\documentclass[11pt,letterpaper]{amsart}

\usepackage{amsmath}
\usepackage{amsthm}
\usepackage{amssymb}
\usepackage{amsfonts}
\usepackage{amsxtra}
\usepackage{fullpage}
\usepackage{amscd}
\usepackage{color} 
\usepackage{bm} 
\usepackage{mathrsfs}
\usepackage{subfigure}
\usepackage{tikz}
\usepackage{tikz-cd}
\usetikzlibrary{cd}
\usepackage{enumerate}
\usepackage[all]{xy}
\usepackage[mathscr]{eucal}
\usepackage{verbatim}

\let\nc\newcommand
\let\renc\renewcommand

\theoremstyle{plain}
\newtheorem{thm}{Theorem}
\newtheorem{prop}[thm]{Proposition}
\newtheorem{cor}[thm]{Corollary}
\newtheorem{lem}[thm]{Lemma}
\newtheorem{conjecture}[thm]{Conjecture}

\theoremstyle{definition}
\newtheorem{defn}[thm]{Definition}

\newtheorem{question}[thm]{Question}

\numberwithin{thm}{section}

\newtheorem*{theorem*}{Theorem}

\nc{\bdm}{\begin{displaymath}}
\nc{\edm}{\end{displaymath}}
\nc{\bthm}{\begin{thm}}
\nc{\ethm}{\end{thm}}
\nc{\blem}{\begin{lem}}
\nc{\elem}{\end{lem}}
\nc{\bcor}{\begin{cor}}
\nc{\ecor}{\end{cor}}
\nc{\bprop}{\begin{prop}}
\nc{\eprop}{\end{prop}}
\nc{\bdef}{\begin{defn}}
\nc{\eddef}{\end{defn}}

\makeatletter
\renewcommand{\subsection}{\@startsection{subsection}{2}{0pt}{-3ex
plus -1ex minus -0.2ex}{-2mm plus -0pt minus
-2pt}{\normalfont\bfseries}} \makeatother

\numberwithin{equation}{section}

\newcommand{\idot}{{\:\raisebox{2pt}{\text{\circle*{1.5}}}}}

\DeclareMathOperator{\Ext}{\mathrm{Ext}}

\DeclareMathOperator{\Ch}{\mathrm{Ch}}

\DeclareMathOperator{\End}{\mathrm{End}}

\DeclareMathOperator{\gr}{\mathrm{gr}}

\DeclareMathOperator{\ann}{\mathtt{Ann}}
\DeclareMathOperator{\Kdim}{\mathrm{Kdim}}
\DeclareMathOperator{\GKdim}{\mathrm{GKdim}}
\DeclareMathOperator{\lKdim}{\mathrm{l-Kdim}}
\DeclareMathOperator{\rKdim}{\mathrm{r-Kdim}}

\newcommand{\beq}{\begin{equation}\label}
\newcommand{\eeq}{\end{equation}}

\newcommand{\iso}{{\;\stackrel{_\sim}{\to}\;}}

\DeclareMathOperator{\Hom}{\mathrm{Hom}}

\nc{\Z}{\mathbb{Z}}
\newcommand{\N}{\mathbb{N}}
\newcommand{\Q}{\mathbb{Q}}
\newcommand{\R}{\mathbb{R}}
\newcommand{\C}{\mathbb{C}}
\newcommand{\h}{\mathfrak{h}}

\nc{\rank}{\textrm{rank} \,}
\nc{\ds}{\dots}

\let\mc\mathcal
\let\mf\mathfrak

\nc{\mbf}{\mathbf}
\nc{\Res}{\mathsf{Res} \, }
\nc{\Ind}{\mathsf{Ind} \, }
\nc{\cont}{\textrm{cont}}

\nc{\msf}{\mathsf}

\nc{\minusone}{-1}
\nc{\minustwo}{-2}
\nc{\Mod}{\mathrm{Mod} \,}
\nc{\ms}{\mathscr}
\nc{\Frac}{\mathrm{Frac} \,}
\nc{\ra}{\rightarrow}
\nc{\hra}{\hookrightarrow}
\nc{\lab}{\label}
\renc{\O}{\mc{O}}
\nc{\Tan}{\mc{T}}
\nc{\ul}{\underline}
\nc{\s}{\mathfrak{S}}
\nc{\g}{\mf{g}}
\nc{\pa}{\partial}
\nc{\tit}{\textit}
\nc{\Maxspec}{\mathrm{Maxspec} \, }
\nc{\gldim}{\mathrm{gl.dim}}
\nc{\rkm}{\mathrm{rk} \, (\mf{m})}
\nc{\sm}{\mathrm{sm}}
\nc{\PD}{\mathbb{PD}}
\nc{\hilb}{\textrm{Hilb}}
\nc{\<}{\langle}
\renc{\>}{\rangle}
\nc{\T}{\mathbb{T}}
\nc{\X}{\mathbb{X}}
\nc{\F}{\mathbb{F}}
\nc{\id}{\msf{id}}
\nc{\A}{\mathbb{A}}
\nc{\Grat}{\mc{Grat}}
\nc{\Squo}[1]{\A^{(#1)}}
\nc{\twist}{\mathrm{twist}}
\nc{\Cd}{\mc{C}}
\nc{\Span}{\mathrm{Span}}
\nc{\Grass}{\mathrm{Gr}}
\nc{\Supp}{\mathrm{Supp}}
\nc{\Irr}{\mathrm{Irr}}
\renc{\o}{\otimes}
\renc{\gr}{\mathsf{gr}}
\nc{\fin}{\mathrm{fin}}
\nc{\aff}{\mathrm{aff}}
\nc{\algD}{\mf{D}}
\nc{\hr}{\mf{h}_{\textrm{reg}}}
\nc{\D}{\mathscr{D}}
\nc{\PIdeg}{\mathrm{PI-degree}}
\nc{\ch}{\mathrm{ch}}
\nc{\ev}{\mathsf{ev}}
\nc{\Stab}{\mathrm{Stab}}
\nc{\Der}{\mathrm{Der}}
\nc{\rightsim}{\stackrel{\sim}{\longrightarrow}}
\nc{\HZ}{H_{\mbf{h},\Z}(\Z_m)}
\nc{\sing}{\mathrm{sing}}
\nc{\dd}{\mathscr{D}}
\nc{\bc}{\mathbf{c}}
\nc{\vc}{\underline{\mathbf{c}}}
\nc{\ba}{\mathbf{a}}
\nc{\reg}{\mathrm{reg}}
\nc{\Amp}{\mathrm{Amp}}
\nc{\Nef}{\mathrm{Nef}}
\nc{\SL}{\mathrm{SL}}
\nc{\Sp}{\mathrm{Sp}}
\nc{\Sym}{\mathrm{Sym}}
\nc{\Mov}{\mathrm{Mov}}
\nc{\Pic}{\mathrm{Pic}}
\nc{\Cs}{\C^{\times}}
\nc{\Nak}[3]{\mf{M}_{{#1}} ({#2},{#3}) }
\nc{\Naka}[2]{\mf{M}({#1},{#2}) }
\nc{\Mtheta}[1]{\mc{M}_{#1}}

\nc{\bw}{\mathbf{w}}

\nc{\bn}{\mathbf{n}}
\nc{\CB}{\mathrm{CB}}
\nc{\GVect}{\Lambda}
\nc{\pZ}{\overline{Z}}
\nc{\Tang}{\mc{T}}
\nc{\K}{\mathbb{K}}

\newcommand{\mr}{\mathrm}

\nc{\red}[1]{\textcolor{red}{#1}}

\begin{document}

\title{Module structure of Weyl algebras}

\author[G. Bellamy]{Gwyn Bellamy}
\address{School of Mathematics and Statistics, University of Glasgow, University Place,
Glasgow, G12 8QQ.}
\email{gwyn.bellamy@glasgow.ac.uk}
%\urladdr{http://www.maths.gla.ac.uk/~gbellamy/}


\begin{abstract}
The seminal paper ``J.T. Stafford, \textit{Module structure of Weyl algebras},  J. London Math. Soc. (2) 18 (1978), no. 3, 429–442" was a major step forward in our understanding of Weyl algebras. Beginning with Serre's Theorem on free summands of projective modules and Bass' Stable Range Theorem in commutative algebra, we attempt to trace the origins of this work and explain how it led to Stafford's construction of non-holonomic simple modules over Weyl algebras. We also describe Bernstein-Lunts' geometric construction of infinite families of non-holonomic simple modules. We recall more recent developments related to Weyl algebras, especially that of parametrizing right ideals in the first Weyl algebra and its relation to Calogero-Moser spaces. Finally, we revisit Stafford's results in the context of quantized symplectic singularities, where they lead naturally to open problems on the behaviour of simple modules.   
\end{abstract}

\maketitle
%\tableofcontents

\section{Introduction}

Weyl algebras were introduced in the early 20th century, in the context of quantum mechanics and differential operators. As the name suggests, the defining relations of these algebras were written down by Hermann Weyl \cite[Equation~(36)]{WeylQuantumArticle} in his foundational article \textit{Gruppentheorie und Quantenmechanik} in 1927, later developed into the book of the same name. The $n$th Weyl algebra $A_n$ is the $\C$-algebra generated by variables $x_1, \ds, x_n, y_1, \ds, y_n$ and satisfying the defining relations
\begin{equation}\label{eq:Weylcommutator}
[x_i,x_j] = [y_i,y_j] = 0, \quad [y_i,x_j] = \delta_{i,j}, \quad \textrm{for all $1 \le i,j \le n$}.
\end{equation}
Here $[A,B]:= AB - BA$ is the commutator bracket. 

In the context of quantum mechanics, letting $x_i$ denote the position operator and $p_i$ the momentum operator, the canonical commutation relation 
\[
[x_i,p_j] = \sqrt{-1} \hbar \delta_{i,j}
\]
underlying Heisenberg's uncertainty principal becomes the relation $[y_j,x_i] = \delta_{i,j}$ if we set $y_i = \sqrt{-1} \hbar^{-1} p_i$. It is in this context that Weyl originally introduced the algebra bearing his name.

Alternatively, one can consider the algebra $\dd(\C^n)$ of differential operators on $\C[x_1, \ds, x_n]$. The ring $\dd(\C^n)$ is generated as a $\C$-algebra by $x_1, \ds, x_n, \partial_1, \ds, \partial_n$, where $\partial_i := \frac{\partial}{\partial x_i}$, and the product rule for differentiation implies that $[\partial_i,x_j] = \delta_{i,j}$. Thus $A_n \iso \dd(\C^n)$ via $y_i \mapsto \partial_i$. 

These different realizations of the Weyl algebra mean that it can be studied from the point of view of $C^*$-algebras, systems of partial differential equations, (complex) symplectic geometry or non-commutative ring theory. The truly interesting mathematics arises from thinking of the Weyl algebra as a bridge between these four disparate fields. This bridge has been extensively exploited and developed in the 98 years since Weyl's article was published. 

This article surveys Stafford's fundamental work in describing the ring-theoretic and module-theoretic properties of Weyl algebras during the decade 1976-1987. In particular, we focus on the paper \textit{Module structure of Weyl algebras} \cite{StaffordModuleStructure} published in the Journal of the London Mathematical Society in 1978 and the subsequent paper \textit{Non-holonomic modules over Weyl algebras and enveloping algebras} \cite{StaffordNonholonomic} published in Inventiones Mathematicae in 1985. These papers beautifully illustrate that Weyl algebras exhibit very surprising (and frankly bizarre) ring-theoretic behaviour confounding intuition coming from commutative ring theory. As we'll see, this behaviour stems largely from the fact that Weyl algebras are \textit{simple rings} - they contain no proper two-sided ideals. 

\subsection*{Outline}

The results in the papers we are surveying are ring-theoretic and their genesis can clearly be seen in the earlier papers \cite{StaffordCompletelyFaithful, StaffordStableStructure} that came out of Stafford's PhD thesis. Since these results provide motivation for the results in \cite{StaffordModuleStructure}, we begin in Section~2 with a summary of these earlier papers. We recall three fundamental results in commutative ring theory due to Serre and Bass. As we will see, Stafford showed that these results have natural counterparts for simple non-commutative rings and are closely related to the basic question of how many elements are required to generate a given module. Section~3 explains how these results were developed and greatly strengthened in the paper \textit{Module structure of Weyl algebras} \cite{StaffordModuleStructure} published in the Journal of the London Mathematical Society. 

As is often the case in mathematics, part of the appeal and beauty of the results in this survey come not just from the fact that they are surprising and deep but that they can be translated into different results in a (seemingly tangentially) related field where they reveal truly surprising and revolutionary consequences. Therefore, after a short Section~4 recalling the definition of holonomic modules over Weyl algebras, we describe in Section~5 Stafford's remarkable results around the existence of non-holonomic simple modules over Weyl algebras. These examples showed that simple modules can behave in very surprising ways and that intuition built from the properties of holonomic modules can be very misleading in general. We also explain how Stafford used these examples to describe the unusual behaviour of simple modules over the enveloping algebra of a semi-simple Lie algebra. Finally, we describe subsequent work of Bernstein-Lunts, where they give a geometric approach to constructing families of examples of non-holonomic simple modules for Weyl algebras. 

In Section~6, we describe more recent developments in the study of Weyl algebras. This section is very subjective and we have decided to concentrate on the remarkable properties of right ideals in the first Weyl algebra. In the final Section~7 we highlight the fact that Stafford's work, viewed in the context of quantized symplectic singularities, raises several interesting open problems. 

\subsection*{Acknowledgements}

 We would like to thank the referee for an exceptionally thorough reading of the survey and many comments that greatly improved the accuracy and presentation. Likewise, our thanks to Ken Brown and Thierry Levasseur for helpful comments on earlier drafts of the survey.

%We would like to thank  for stimulating discussions about . 

\section{Prehistory - The stable Range Theorems}\label{sec:prehistory1}

Our story begins not with the Weyl algebra but rather at a key moment in the establishment of modern commutative algebra. In the late 1950's, J. P. Serre did seminal work in translating concepts in topology into the language and setting of algebra. A fundamental observation of Serre \cite{SerreProjective} along these lines is that the category of (algebraic) vector bundles on an affine variety is equivalent to that of finitely generated projective modules over the corresponding commutative Noetherian ring. He goes on to prove in the algebraic setting several results that were known previously in the geometric setting. The most important of these is \textit{Serre's Theorem} saying that projective modules have a free summand if they are sufficiently large.

\begin{thm}[Serre's Theorem]\label{thm:Serre1} 
	Let $R$ be a commutative Noetherian domain with $\Kdim R = n$ and $P$ a finitely generated projective module of rank $> n$. Then $P \cong P' \oplus R$ for some (projective) submodule $P'$. 
\end{thm}

Here $\Kdim R$ denotes the Krull dimension of the ring $R$ \cite[Chapter~5]{MatCom} and the rank of a finitely generated module $M$ is defined to be $\dim_K (K \o_R M)$, where $K$ is the field of fractions of the domain $R$. The above result is \cite[Th\'eor\`eme~1]{SerreProjective}.  As noted at the end of Serre's paper, this result has applications to understanding the structure of the Grothendieck group of projective $R$-modules. Later Bass \cite{BassStable} undertook a far more systematic study of the Grothendieck group $K_0$, together with the higher $K$-group $K_1$. In this context, Bass establishes two other key results closely related to Serre's Theorem. The first of these is \cite[Theorem~11.1]{BassStable}: 

\begin{thm}[Stable Range Theorem]
	Let $R$ be a commutative Noetherian ring with $\Kdim R = n$, and suppose $R = \sum_{i = 1}^{n+2} a_i R$. Then there exist $f_i \in R$ such that 
	\[
	R = \sum_{i = 1}^{n+1} (a_i + a_{n+2} f_i) R. 
	\] 
\end{thm}

Bass was motivated to prove this theorem in order to understand the properties of $GL(R,n)$ as $n \to \infty$. The terminology ``stable range" refers to $n \in \N$ for which properties of $GL(R,m)$ stabilize (i.e., are independent of $m$) for all $m \ge n$, by analogy with stable homotopy theory. 

\begin{thm}[Cancellation Theorem]
	Let $R$ be a commutative Noetherian domain with $\Kdim R = n$ and $M$ a projective module of rank $> n$. If $Q$ is another finitely generated projective module and $N$ any $R$-module such that $M \oplus Q \cong N \oplus Q$ then $M \cong N$. 
\end{thm}

The cancellation theorem also appears to be due to Bass from his work on the stable range for $GL(R)$, see \cite[Theorem~9.3]{BassStable}. 

%Though not immediately apparent, both Serre's Theorem (which appears as Theorem~8.2 in Bass' paper) and the Cancellation Theorem can be viewed as consequences (CHECK!) of the Stable Range Theorem

One would imagine that all this is rather far removed from the paper \cite{StaffordModuleStructure} by Stafford on modules over Weyl algebras. However, a glance at the papers \cite{StaffordCompletelyFaithful, StaffordStableStructure} resulting from his PhD thesis make the connection very clear. In these papers, Stafford explores the extent to which the above three results make sense for non-commutative simple (left or right) Noetherian rings. As we'll see, much stronger versions of these results can be shown to hold when the simplicity hypothesis is assumed. It is surprising to the author that one would even consider these fundamental results in commutative algebra in the context of simple non-commutative rings given that the two classes are in a sense diametrically opposite. 

As a precursor to the Stable Range Theorem, Stafford shows that: 

\begin{thm}\cite[Theorem~1.3]{StaffordCompletelyFaithful}\label{thm:Thm2.4}
	Suppose $n < \infty$ and $R$ is a ring with\footnote{Following \cite{StaffordStableStructure}, $\lKdim R$ denotes the Krull dimension of $R$ considered as a left $R$-module, and similarly for $\rKdim R$.} $\rKdim R \ge n$. If $M$ is a completely faithful Noetherian right $R$-module such that $\Kdim M = n-1$ then $M$ can be generated by $n$ elements. 
\end{thm}

Here $\Kdim M$ denotes Krull dimension for non-commutative rings in the sense of Rentschler-Gabriel, see e.g., \cite[Chapter~6]{MR}. This is a measure of the size of the poset of submodules of $M$ and $\Kdim R$ agrees with classical Krull dimension for commutative Noetherian rings $R$. A module $M$ is faithful if its annihilator $\{ r \in R \, | \, rm = 0 \, \forall m \in M \}$ equals zero and $M$ is completely faithful if all non-zero submodules of all quotients of $M$ are faithful. Notice that if $R$ is simple then every non-zero module is completely faithful.  

\begin{cor}\cite[Corollary~1.5]{StaffordCompletelyFaithful}\label{cor:1.5}
	If $R$ is a simple right Noetherian ring with Krull dimension $n$ then any right ideal $I$ of $R$ can be generated by $n+1$ elements. 
\end{cor}

In particular, these results hold for the Weyl algebra $A_n$, though we will soon see that much more can be said in this case. In the subsequent paper \cite{StaffordStableStructure} Stafford goes on to establish a full-blown Stable Range Theorem for modules over a simple Noetherian ring, strengthening Theorem~\ref{thm:Thm2.4}. 

\begin{thm}\cite[Theorem~2.2]{StaffordStableStructure}\label{thm:simplestable}
Let $R$ be a right Noetherian simple ring with $\rKdim R \ge n$. Let $M$ be a finitely generated right $R$-module with $\Kdim M = n-1$ and suppose that $M = \sum_{i = 1}^{m+1} a_i R$ with $m \ge n$. Then there exist $f_i \in R$ such that 
\[
M = \sum_{i = 1}^{m} (a_i + a_{m+1} f_i) R. 
\] 
\end{thm}

From this theorem, Stafford goes on to prove versions of the Stable Range Theorem for right ideals in $R$ \cite[Theorem~2.4]{StaffordStableStructure} including for the ring itself. He then proves versions of Serre's Theorem and the Cancellation Theorem for simple rings. Both proofs rely heavily on the Stable Range Theorem~\ref{thm:simplestable}. 

\begin{thm}\cite[Theorem~4.3]{StaffordStableStructure}
		Let $R$ be a simple left Noetherian ring with $\lKdim R = n$ and $M$ a torsion-free right module of rank $\ge n+2$. Then $M \cong M' \oplus R$ for some submodule $M'$. 
\end{thm}

The fact that $R$ is assumed to be left Noetherian and $M$ a right module is \textit{not} a mistake in the above theorem. Here the rank of a torsion free module $M$ is defined to be the smallest integer $r$ such that $M$ embeds into $R^{(r)}$. The Cancellation Theorem in this context says that:

\begin{thm}\cite[Theorem~4.5]{StaffordStableStructure}
	Let $R$ be a simple left Noetherian ring with $\lKdim R = n$ and $N$ a right $R$-module which has a torsion-free direct summand of rank $\ge n+2$. Let $P$ be a finitely generated projective $R$-module and $N'$ any $R$-module such that $N \oplus P \cong N' \oplus P$. Then $N \cong N'$. 
\end{thm}

%Stafford goes on to recombine these reults with the original commutative situation of Serre and Bass by considering projective modules over $R[x_1, \ds, x_m]$ for $R$ a simple ring. 

%\begin{thm}\cite[Theorem~5.4]{StaffordStableStructure}
%Let $R$ be a simple Noetherian ring with uniform dimension $y$ and $\lKdim R = n$. Let $S = R[x_1, \ds, x_m]$ and $P$ be a finitely generated projective right $S$-module with $\rk(P) \ge \max (n+2,m/y + 1)$. Then $P \cong P' \oplus S$. 
%\end{thm}

%\begin{thm}\cite[Theorem~5.6]{StaffordStableStructure}
%Let $R$ be a simple Noetherian ring with uniform dimension $y$ and $\lKdim R = n$. Let $S = R[x_1, \ds, x_m]$ and $P$ be a finitely generated projective right $S$-module with $\rk(P) \ge 1 + \max (n+2,m/y + 1)$. If $P'$ and $P''$ are finitely generated projective $S$-modules such that $P \oplus P'' \cong P' \oplus P''$ then $P \cong P'$. 
%\end{thm}

In the case of the Weyl algebra $A_n$, which is treated in Section~6 of \cite{StaffordStableStructure}, the above results are hugely strengthened. This is the first point where we see the exceptional properties of the Weyl algebra emerging and the results diverge from the uniform ones stated previously. Since these results are all superseded by those in \cite{StaffordModuleStructure}, we do not recall them here except to say that at this point Stafford is able to show that each right ideal of $A_n$ can be generated by at most $5$ elements. In particular, this bound does not depend on the Krull dimension of $A_n$; compare with Corollary~\ref{cor:1.5}.  

%\begin{thm}\cite[Theorem~6.1]{StaffordStableStructure}
%	Let $M$ be a finitely generated torsion-free $A_n$-module with $\rk(M) \ge 4$. Then $M \cong M' \oplus A_n$. 
%\end{thm}

%\begin{thm}\cite[Theorem~6.3]{StaffordStableStructure}
%	Let $N$ be a right $A_n$-module such that $N = N' \oplus M$ where $M$ is a torsion-free finitely generated right $A_n$-module with $\rk(M) \ge 5$. Suppose $L$ and $P$ are $A_n$-modules with $P$ finitely generated, projective, and $N \oplus P \cong L \oplus P$. Then $N \cong L$.  
%\end{thm}

%In particular the above result implies that if $P$ is a projective finitely generated $A_n$-module with $\rk(P) \ge 5$ then $P$ is free \cite[Corollary~6.4]{StaffordStableStructure}. 

%Finally, we get to first iteration of the Stable Range Theorem for Weyl algebras. 

%\begin{thm}\cite[Theorem~6.5]{StaffordStableStructure}
%	Let $I$ be a left ideal of $A_n$ and suppose $I = \sum_{i = 1}^{r+1} A_n a_i$ with $r \ge 5$. Then there exist $f_i \in A_n$ such that $I = \sum_{i = 1}^{r} A_n (a_i + f_i a_{r+1})$. 
%\end{thm}

%In particular, this implies that any left or right ideal in $A_n$ can be generated by just $5$ elements! This is the first hint at the surprising results to come in \cite{StaffordModuleStructure}.

\section{Module structure of Weyl algebras}

Krull dimension plays a crucial role in the proof of the results described in Section~\ref{sec:prehistory1}. This is evident even from the statement of the results, which is encoded in terms of Krull dimension. However, Krull dimension does not appear in the statement of the results for the Weyl algebra (which is partly why they are so surprising) and neither is it used in the proof. 

The key technical result is the following version of the Stable Range Theorem. 

\begin{thm}\cite[Theorem~3.1]{StaffordModuleStructure}\label{thm:SRTWeylfinal}
	Let $I = a A_n + b A_n + c A_n$ be a right ideal of $A_n$ and let $d_1, d_2$ be arbitrary non-zero elements of $A_n$. Then there exist $f$ and $g$ in $A_n$ such that 
	\[
	I = (a + cfd_1) A_n + (b + cg d_2) A_n. 
	\] 
\end{thm}

In particular, this implies that every right ideal in $A_n$ can be generated by just two (!!) elements. One should compare this with Theorem~\ref{thm:simplestable}, noting that the Krull dimension of $A_n$ is $n$ \cite[Theorem~6.6.15]{MR}. In terms of PDEs with polynomial coefficients, Theorem~\ref{thm:SRTWeylfinal} implies that if $L_1, \ds, L_k \in A_n$ then there exist $P,Q \in A_n$ such that a function $f$ is a solution to the system of equations $L_1 f = \cdots = L_k f = 0$ if and only if $P f = Q f = 0$. Theorem~\ref{thm:SRTWeylfinal} is an abstract existence result. However, A. Leykin \cite{LeykinStafford} has developed an algorithm implemented in Macaulay2 \cite{Macaulay2} to explicitly compute two generators of a given right ideal. The following analogues of Serre's Theorem and the Cancellation Theorem are deduced, more or less directly, from Theorem~\ref{thm:SRTWeylfinal}. 

\begin{thm}\cite[Theorem~3.3]{StaffordModuleStructure}
	Let $M$ be a finitely generated left $A_n$-module. Then $M \cong M' \oplus A_n^{(s)}$, where $M'$ is a module of rank $\le 1$. If $M$ is torsion-free then $M'$ is isomorphic to a left ideal of $A_n$. 
\end{thm}

\begin{thm}\cite[Theorem~3.6]{StaffordModuleStructure}
	(a) Let $M$ be a finitely generated left $A_n$-module with rank $M \ge 2$. Suppose that $M \oplus A_n \cong N \oplus A_n$ for some module $N$. Then $M \cong N$. 
	
	(b) Let $P$ be a projective left $A_n$-module. Then either $P$ is free or $P \cong I$, a left ideal of $A_n$.  
\end{thm}

Notice that if we take $I = A_n$ in Theorem~\ref{thm:SRTWeylfinal} then it says that for arbitrary non-zero elements $d_1, d_2 \in A_n$ there exist $f,g \in A_n$ such that $A_n = f d_1 A_n + g d_2 A_n$. Towards the end of the paper Stafford returns to the question of how many elements are required to generate a module. If $M$ is a finitely generated torsion free module of rank $r$ then \cite[Theorem~3.9]{StaffordModuleStructure} says that either $M \cong A_n^{(r)}$ or $M$ can be generated by $r+1$ elements and no fewer. For a finitely generated torsion module $T$, it is shown in \cite[Theorem~3.7]{StaffordModuleStructure}  that $T$ can be generated by two elements. Stafford conjectures:  

\begin{conjecture}\label{conj:torsion}
Every finitely generated torsion $A_n$-module is cyclic. 
\end{conjecture}

%Let $d \in A_n$ be non-zero. Then there exists $f \in A_n$ such that $A_n = d A_n + fd A_n$. 

This conjecture is stated in a different (and slightly stronger) form as \cite[Conjecture~3.8]{StaffordModuleStructure} in the paper. As far as the author is aware, this conjecture is still open. 

Stafford's proof of Theorem~\ref{thm:SRTWeylfinal} is by induction on $n$, beginning with the first Weyl algebra. He makes clever use of partial quotient rings of $A_n$ in order to lift results from $A_r$ for $r < n$. The arguments are elementary in the sense that he works directly with elements of these partial quotient rings and does not required any results from his earlier papers. 

\section{Holonomic $\dd$-modules}\label{sec:prehistory2}

We have seen that there is a clear path from commutative algebra and $K$-theory, through analogous results for simple Noetherian rings, to the generation and Stable Range theorems for Weyl algebras. An important aspect of the representation theory of Weyl algebras which has not yet made an appearance is the notion of holonomic modules. Though the notion appears to have its origin in the study of ``maximally overdetermined systems of partial differential equations" (see \cite{SKK}), it is more clearly motivated by asking the basic question: what is the ``size" of simple modules over the Weyl algebra? Before making this precise, consider the following elementary fact. If $M$ is a non-zero $A_n$-module then $\dim_{\C} M = \infty$. Indeed, suppose that $\dim_{\C} M < \infty$. Then $x_1,y_1$ act as linear transformations $X_1,Y_1$ of $M$ such that $[Y_1,X_1] = \mathrm{Id}_M$. Taking trace of both sides of this equation gives $0 = \dim_{\C} M$, a contradiction. Therefore one needs a way to measure size beyond naive dimension as a $\C$-vector space. 

The correct notion of size is Gelfand-Kirillov dimension. We will not recall the definition here. Suffice it to say, this is a measure of the rate of growth of finitely generated modules over an affine (i.e., finitely generated) algebra. A key result is Bernstein's inequality \cite[Theorem~1.3]{BernsteinContinuation}, which states that if $M$ is a non-zero finitely generated module over the Weyl algebra $A_n$ then 
\begin{equation}\label{eq:Bernsteinineq}
2 \GKdim (M) \ge \GKdim (A_n) = 2n \quad \textrm{i.e.,} \quad \GKdim (M) \ge n.
\end{equation}
The holonomic modules are precisely those of smallest Gelfand-Kirillov dimension i.e., the modules $M$ with $\GKdim (M) = n$. In many ways, the (abelian) category of holonomic modules behaves like the category of finite-dimensional modules over a finite-dimensional algebra. For instance, holonomic modules have finite length and $\Hom_{A_n}(M,N)$ is a finite-dimensional $\C$-vector space when $M,N$ are holonomic. Moreover, it has been shown that every holonomic module is cyclic, which to the author is a very surpring result; see e.g.,  \cite[Corollary~10.2.6]{PrimerDmodules}. This is relevant for Conjecture~\ref{conj:torsion} since holonomic modules are torsion. 

The standard example of a holonomic module is $\C[x_1, \ds, x_n]$, with $A_n$ acting as differential operators. It is easy to see that this module is both holonomic and simple. The ring itself gives an example of a module which is not holonomic. There is an abstract classification of simple holonomic modules for $A_n$ in terms of minimal extensions of integrable connections; see e.g., \cite[\S3.4]{HTT}.

\section{Existence of non-holonomic simple modules}

Though we have seen an example of a simple holonomic module and they are classifiable in a certain sense, it is much harder to write down a simple module that is not holonomic. In fact, the expectation in the 1970s was that all simple modules over $A_n$ should be holonomic. In the literature, one can see this as a special case of a more general expectation of a relation between the Krull dimension and the GK-dimension of a module. Precise statements can be found in Problem~1 (Chapter~1,~\S9) of the book \cite{BjorkBook} by J.-E. Bj\"ork and Section~7 of \cite{MCWeyl}. 

Therefore, it was very suprising when Stafford showed: 

\begin{thm}\label{thm:nonholonomic}\cite[Theorem~1.1]{StaffordNonholonomic}
There exists a cyclic, maximal right ideal $a A_n$ in $A_n$. In particular, the simple $A_n$-module $N = A_n/ a A_n$ has Gelfand-Kirillov dimension $2n-1$. 
\end{thm}

If you want to check that this contradicts the above mentioned expectations, note that a simple module has Krull dimension zero. Stafford goes on to show that these non-holonomic modules behave very differently to holonomic modules. 

To illustrate this, let us first recall the conjecture by I. M. Gelfand that motived J. Bernstein to introduce, and study, holonomic modules over the Weyl algebra. Let $P$ be a polynomial in $n$ variables with real coefficients and $\lambda \in \C$ with $\mathrm{Re} \, \lambda > 0$. If $\Theta$ is an open subset of $\R^n$ on which $P$ is positive then we can define the (generalized) function $P^{\lambda} \colon \R^n \to \C$ by $P^{\lambda}(x) = P(x)^{\lambda}$ for $x \in \Theta$, and zero outside $\Theta$. In the International Congress of Mathematicians, Amsterdam, 1954, Gelfand conjectured in his lecture \cite{GelfandICM} that $P^{\lambda}$ extends to a function on $\R^n \times \C$ that is meromorphic on $\C$. In fact, he made a strong conjecture on the meromorphic extension of certain distributions associated to $P$. A proof of this (more general) conjecture was given by S. I. Gelfand--J. Bernstein \cite{BernsteinGelfand} and M.~Atiayh \cite{AtiyahGelfand}, both relying on Hironaka's Theorem on the existence of resolutions of singularities. 

In 1972, Bernstein used holonomic modules to give a pure algebraic proof \cite[Theorem~1]{BernsteinContinuation} of the conjecture, as stated above. This was done by proving the existence \cite[Theorem~1']{BernsteinContinuation} of the ``Bernstein-Sato'' polynomial associated to $P$. The key step in the proof of the latter result is to show that if $M$ is holonomic and $f \in \C[x_1, \ds, x_n]$ non-zero then the localization $M[f^{-1}]$ is still holonomic. In particular, it is finitely generated. This is a central result in the theory of holonomic $\dd$-modules since it allows one to prove that holonomicity is preserved under (derived) pull-back and push-forward for arbitrary morphisms. See \cite[Chapter~8]{KrauseLenagan} for a detailed exposition of Bernstein's proof of Gelfand's conjecture. 

Contrast Bernstein's result on the localization of holonomic modules with the following corollary by Stafford of Theorem~\ref{thm:nonholonomic}.

\begin{cor}\cite[Corollary~1.2]{StaffordNonholonomic}\label{cor:simplelocalization}
	Let $\alpha = x_1 + y_1 x_2 y_2 + x_2 + y_2 \in A_2$ and $M = A_2 / \alpha A_2$. Then $M$ is a simple $A_2$-module such that $M \o_{\C[x_1,x_2]} \C[x_1^{\pm 1},x_2]$ is not finitely generated as a right $A_2$-module. 
\end{cor}

Another result already mentioned is that the space of homomorphisms between holonomic modules is finite-dimensional over $\C$. This holds more generally for all Ext-groups between holonomic modules. Therefore, if we take $\alpha = x_1 + y_1 x_2 y_2 + x_2 + y_2 \in A_2$ and $M := A_2 / \alpha A_2$ as in Corollary~\ref{cor:simplelocalization}, it is very surprising that:

\begin{cor}\cite[Corollary~1.3]{StaffordNonholonomic}
 $\Ext^1_{A_2}(M,M)$ is an infinite-dimensional $\C$-vector space. 
\end{cor}

Perhaps even more strangely, V. Bavula gives in \cite[Theorem~2.2]{BavulaQuestions} a simple holonomic module $S$ such that $\Ext^1_{A_n}(N,S)$ is infinite-dimensional as a complex vector space. Here $N = A_n/ a A_n$ is the simple module introduced in Theorem~\ref{thm:nonholonomic}. This implies that there exist (uncountably many) indecomposable modules $T$ of length two with submodule isomorphic to $S$ and quotient $N$. This answered positively a question raise by T. Lenagan \cite{LenaganQuestion}. In the opposite direction, G. Perets \cite{Perets} gave an example of an indecomposable module of length two with non-holonomic simple submodule $N$ and simple holonomic quotient.

\subsection{Applications to enveloping algebras}\label{sec:applicationstoenvelopin}

Bernstein's inequality is also known to hold \cite[Theorem~9.11]{KrauseLenagan} for the enveloping algebra of a finite dimensional (complex) Lie algebra $\mf{g}$; a result due to O. Gabber. However, the statement of the inequality must be modified to take into account the fact that the enveloping algebra is not a simple ring. Namely, Gabber showed that if $M$ is a finitely generated $U(\mf{g})$-module with annihilator $\ann_U(M) \lhd U(\mf{g})$ then:
\begin{equation}\label{eq:BernsteinLie}
2 \GKdim(M) \ge \GKdim(U(\mf{g}) / \ann_U(M)). 
\end{equation}
A simple $U(\mf{g})$-module $M$ is \textit{holonomic} if we have equality in \eqref{eq:BernsteinLie}. More generally, a module $M$ is holonomic if it has finite length and all simple subquotients are holonomic. 

If $\mf{g}$ is a semi-simple Lie algebra then the Conze embedding gives an embedding of each primitive central quotient $U(\mf{g})/P$ of $U(\mf{g})$ into some Weyl algebra $A_n$. Stafford considers the Conze embedding of a primitive factor $U(\mf{sl}_2 \times \mf{sl}_2) / P$ of $U(\mf{sl}_2 \times \mf{sl}_2)$ into $A_2$. He shows that one can choose $\beta \in A_2$ (unfortunately one cannot use the element $\alpha$ appearing in Corollary~\ref{cor:simplelocalization} for this) such that $A_2 / \beta A_2$ is simple and $U(\mf{sl}_2 \times \mf{sl}_2) / (\beta A_2 \cap U(\mf{sl}_2 \times \mf{sl}_2))$ is a simple $U(\mf{sl}_2 \times \mf{sl}_2)/P$-module. 

\begin{thm}\label{thm:enveloping1}\cite[Theorem~2.7]{StaffordNonholonomic}
	There exists $\beta \in A_2$ such that $U(\mf{sl}_2 \times \mf{sl}_2) / (\beta A_2 \cap U(\mf{sl}_2 \times \mf{sl}_2))$ is a non-holonomic simple $U(\mf{sl}_2 \times \mf{sl}_2)$-module. 
\end{thm}

Since $U(\mf{g})$ is a Hopf algebra, the tensor product $E \o_{\C} M$ of any two $U(\mf{g})$-modules $E$ and $M$ is again a $U(\mf{g})$-module. The following is a truly remarkable result, which beautifully illustrates the strange behaviour of non-holonomic simple modules. 

\begin{thm}\label{thm:enveloping2}\cite[Theorem~4.1]{StaffordNonholonomic}
	There exist simple right $U(\mf{sl}_2 \times \mf{sl}_2)$-modules $E$ and $M$, with $E$ finite-dimensional, such that $E \otimes_{\C} M$ has infinite length. 
\end{thm}

Given the importance of projective functors \cite{BGTensor} in the representation theory of simple Lie algebras, Theorem~\ref{thm:enveloping2} has motivated a great deal of research into finding classes of simple modules $M$ for which  $E \otimes M$ has finite length for all finite-dimensional representations $E$; see e.g., \cite[\S1.3]{ProjFunctors}. Stafford's approach is very explicit, involving concrete computations in $U(\mf{sl}_2 \times \mf{sl}_2)$. It raises the question of whether Theorems~\ref{thm:enveloping1} and \ref{thm:enveloping2} hold for any semi-simple Lie algebra. S.~C.~Coutinho shows in \cite{CoutinhoNonholosemisimple} that this is indeed the case, except possibly when $\mf{g} = \mf{sl}_2$. 

\subsection{The generalization by Bernstein-Lunts}

A couple of years after Stafford's work \cite{StaffordNonholonomic} showing the existence of non-holonomic modules for the Weyl algebra, J. Bernstein and V. Lunts \cite{BLnonholo} used a very different, geometric, method to construct non-holonomic simple modules. They showed that any sufficiently generic $d \in A_n$ gives rise to a simple module $A_n / d A_n$. In this section we describe more precisely their result. 

The link to geometry is made via the characteristic variety of an $A_n$-module, allowing one to define holonomic modules via geometry. This approach requires some preparation but provides in return finer information about the module. The Weyl algebra can be endowed with the \textit{Bernstein filtration} $A_n = \bigcup_{i \ge 0} \mc{B}_i A_n$, where $\mc{B}_i A_n$ is the $\C$-linear span of all polynomials of degree at most $i$ in the generators $x_j,y_k$. Thus, $\mc{B}_0 A_n = \C$, $\mc{B}_1 A_n = \mathrm{span}_{\C} \{ 1, x_1, \ds, x_n, y_1, \ds, y_n \}$ and $(\mc{B}_i A_n ) \cdot (\mc{B}_j A_n)  = \mc{B}_{i+j} A_n$. The associated graded algebra is commutative:
\begin{equation}\label{eq:Bernsteinassgr}
P := \C[x_1, \ds, x_n, \xi_1, \ds, \xi_n] \iso \gr_{\mc{B}} A_n := \bigoplus_{i \ge 0} \mc{B}_{i} A_n / \mc{B}_{i-1} A_n,
\end{equation}
where $\mc{B}_{-1} A_n := \{ 0 \}$ and $\xi_i$ maps to the class of $y_i$ in $\mc{B}_{1} A_n / \mc{B}_{0} A_n$. This isomorphism allows one to use the commutator bracket $[ - , - ]$ on $A_n$ to define a Poisson bracket $\{ - , - \}$ on the polynomial ring $P$. The commutation relations \eqref{eq:Weylcommutator} imply that 
\[
\{ \xi_i,\xi_j \} = \{ x_i, x_j \} = 0, \quad \{  \xi_i,x_j \} = \delta_{i,j},  \quad \textrm{for all $1 \le i,j \le n$},
\]
so that, up to sign, one can view $P$ with this Poisson bracket as the algebra of functions on the symplectic vector space $\C^{2n}$. 

Let $M$ be a finitely generated left $A_n$-module. We say that a filtration $\mc{B}_{\idot} M$ on $M$ is a \textit{good filtration} if 
\[
(\mc{B}_i A_n )(\mc{B}_j M) \subset \mc{B}_{i+j} M, \quad \mc{B}_i M = \{ 0 \} \textrm{ for $i \ll 0$}, \quad M = \bigcup_{i \in \Z} \mc{B}_i M,
\]
and, most importantly, 
\[
\gr_{\mc{B}} M := \bigoplus_{i \in \Z} \mc{B}_{i} M / \mc{B}_{i-1} M
\]
is finitely generated over $P$. One can easily show that any finitely generated $A_n$-module can be equipped with a good filtration, though good filtrations are far from unique. The \textit{characteristic variety} $\mathrm{Ch}_{\mc{B}}(M) \subset \C^{2n}$ of $M$ is the support of the finitely generated $P$-module $\gr_{\mc{B}} M$. As a (reduced) subvariety of $\C^{2n}$, $\mathrm{Ch}_{\mc{B}}(M)$ is independent of the choice of good filtration. We recall that an ideal $I$ in the Poisson algebra $P$ is \textit{involutive} if $\{ I , I \} \subset I$. The following deep result is central in the theory of $\dd$-modules.  

\begin{thm}[Gabber's Theorem \cite{Gabber}]\label{thm:Gabber}
	The ideal $\sqrt{\ann \gr_{\mc{B}} M}$ is involutive.  
\end{thm}

The importance of Gabber's Theorem is that it implies that $\Ch_{\mc{B}}(M) \subset \C^{2n}$ is a coistropic subvariety. This implies that $2 \dim \mathrm{Ch}_{\mc{B}}(M) \ge \dim \C^{2n}$. Since it is easily shown that $\dim \mathrm{Ch}_{\mc{B}}(M) = \GKdim (M)$, we see that Bernstein's inequality \eqref{eq:Bernsteinineq} can deduced from Gabber's Theorem~\ref{thm:Gabber}. %be viewed as a consequence of 

The isomorphism \eqref{eq:Bernsteinassgr} is graded, where $P$ is graded with $\deg x_i = \deg \xi_j = 1$.  We let $P_k$ denote the finite-dimensional subspace of homogeneous polynomials of degree $k$.  

\begin{thm}\label{thm:BL1}\cite[Theorem~A]{BLnonholo} 
 For $n = 2$ and $k \ge 4$, if $d \in \mc{B}_k A_n$ is any element whose symbol $s = \sigma(d) \in P_k$ is generic then $A_n / d A_n$ is irreducible. 
\end{thm}

As in \cite[Remark~1]{BLnonholo}, here ``generic'' means that the set $N \subset P_k$ of polynomials $s$ for which the conclusion of Theorem~\ref{thm:BL1} holds satisfies the following condition: Its complement $P_k \setminus N$ can be covered by a countable number of hypersurfaces. 

Using the fact that $\mathrm{Ch}_{\mc{B}}(M)$ is a homogeneous coisotropic subvariety of $\C^{2n}$, the above theorem is easily deduced from the corresponding geometric statement outlined in Theorem~\ref{thm:BL2} below. Following Bernstein-Lunts, we say that a subvariety $Y \subset \C^{2n}$ is a \textit{minimal homogeneous coistropic subvariety} if it is a homogeneous coistropic subvariety that does not properly contain any other (closed) homogeneous coistropic subvariety (they use the term ``involutive'' rather than coistropic). Using the non-existence of algebraic solutions to a certain generic first order partial differential equation \cite[Theorem~1]{BLnonholo}, they show that:

\begin{thm}\label{thm:BL2}\cite[Theorem~A']{BLnonholo}
 For $n = 2$ and $k \ge 4$ and generic $s \in P_k$, the irreducible subvariety $V(s) \subset \C^{2n}$ is a minimal homogeneous coistropic subvariety. 
 \end{thm}

Again, ``generic'' means that the set $N \subset P_k$ of polynomials $s$ for which the conclusion of Theorem~\ref{thm:BL2} holds has the property that its complement $P_k \setminus N$ can be covered by a countable number of hypersurfaces. Bernstein-Lunts conjecture that the statement of Theorem~\ref{thm:BL2}, and hence of Theorem~\ref{thm:BL1} too, should hold for all $n \ge 2$ and $k \ge 3$. This conjecture was soon confirmed by V. Lunts \cite{LuntsFoliation} for $n \ge 2, k \ge 4$. The case $n  = 2$ and $k = 3$ was done by T.~McCune \cite{McCune} and the induction argument of \cite{LuntsFoliation} applied to this case also deals with $n \ge 2$ and $k = 3$. 
 
Since $A_n$ is a ring of differential operators it can be given the filtration $\mc{F}_{\idot} A_n$ by order of differential operators. Here $\mc{F}_i A_n$ is the subspace spanned by all differential operators of order at most $i$ so that $\mc{F}_0 A_n = \C[x_1, \ds, x_n]$ and $y_i \in \mc{F}_1 A_n$. One can ask what happens if we replaced the Bernstein filtration by the order filtration. The associated graded ring is again the polynomial ring $P$, though one should now think of $P$ as the ring of functions on the cotangent bundle $T^* \C^n$ (the Poisson structure is the same, but the grading is different). Again, every finitely generated module $M$ has a characteristic variety $\Ch_{\mc{F}}(M) \subset T^* \C^n$, coistropic by Gabber's Theorem~\ref{thm:Gabber}. We note that under the obvious identification $T^* \C^n = \C^{2n}$, $\Ch_\mc{B}(M) \neq \Ch_{\mc{F}}(M)$ though they have the same dimension. 

Bernstein-Lunts consider polynomial vector fields $\xi = \sum_{i = 1}^n r_i(x) \partial_i$ and say that a pair $(\xi,f)$, where $f \in \C[x_1, \ds, x_n]$, is very generic if it satisfies certain technical genericity conditions (these are conditions (*) and (**) of \cite[Proposition~6]{BLnonholo}). If $(\xi,f)$ is very generic then they show that $A_n / (f + \xi) A_n$ is a simple module. The difficulty is to prove that there exist any very generic pairs $(\xi,f)$. They show that this is the case when $n = 2$. 
 
 \begin{thm}\label{thm:BL3}\cite[\S~4.3]{BLnonholo}
 If $n = 2$ then there exist very generic pairs $(\xi,f)$, where the coefficients $r_i(x)$ of $\xi$ have degree $\ge 2$. 
 \end{thm}
 
 Again, they conjecture that the statement of Theorem~\ref{thm:BL2} holds for all $n \ge 2$. This conjecture was confirmed by  J. Vit\'orio Pereira \cite{PereiraChar}. 
 
\section{Subsequent developments}\label{sec:subsequent}

As noted in the introduction, one can view the Weyl algebra as the ring of differential operators on $\C^n$. From this point of view, the Weyl algebra is just a special case of rings $\dd(X)$ of differential operators on smooth affine varieties $X$. The study of modules over these rings (or rather modules over the corresponding sheaf of differential operators in the case of non-affine varieties) is a vast area - this is the theory of $\dd$-modules. We won't touch on this subject here. See, for instance, the books \cite{PrimerDmodules,HTT} for excellent introductions. However, one can ask to what extent Stafford's results extend to the rings $\dd(X)$. This question was considered by S.~C.~Coutinho and M.~Holland\footnote{See Section~4 of \cite{CoutinhoHollandModule} for a number of open conjectures in this direction.} \cite{CoutinhoHollandModule}. In particular, they showed that every right ideal in $\dd(X)$ can be generated by $3$ elements; it is a question of J.-E. Bj\"ork \cite[Problem~2.10]{BjorkBernstein} whether $2$ elements suffice. 

Somewhat tangential to the general theory of $\dd$-modules, Stafford considers in \cite{StaffordEdnomorphisms} one-sided ideals in the first Weyl algebra. Since we are working over the complex numbers, the Weyl algebra $A_1$ is hereditary, i.e., has (left and right) global dimension one, and hence every one-sided ideal is projective. If $0 \neq I \subset A_1$ is a right ideal then the fact that $A_1$ is simple and hereditary implies that $A_1$ is Morita equivalent to $E = \End_{A_1}(I)$; it is a generator because $I^* I$ is a non-zero two-sided ideal of $A_1$ and hence $I^* I = A_1$. Stafford shows \cite[Theorem~A]{StaffordEdnomorphisms} that $E \cong A_1$ if and only if $I$ is a cyclic right ideal. As a consequence:

\begin{cor}\cite[Corollary~B]{StaffordEdnomorphisms}\label{cor:Moritaequiv}
Let $P$ and $Q$ be right ideals of $A_1$. Then $\End_{A_1}(P) \cong \End_{A_1}(Q)$ if and only if $P = t \sigma(Q)$ for some $t \in D(A_1)$, the division ring of fractions of $A_1$, and $\sigma \in \mr{Aut}_{\C}(A_1)$, the group of $\C$-algebra automorphisms of $A_1$. 
\end{cor}

This remarkable result raises more questions than it answers. How many isomorphism classes of right ideals are there in the first Weyl algebra? What is $\mr{Aut}_{\C}(A_n)$?... 

The first question has led to an explosion of work connecting the first Weyl algebra to integrable systems, rational Cherednik algebras and non-commutative projective geometry. We only touch on the beginnings of this story. Two right ideals $P,Q$ of $A_1$ are said to be isomorphic if there exist non-zero $s,t \in A_1$ such that $s P = t Q$. Since $A_1$ is a Noetherian domain, this is the same as saying that they are isomorphic as $A_1$-modules. R.~Cummings and M.~Holland \cite{CanningsHolland} described a certain (infinite) set $\mc{C}$ that parametrizes explicitly the isomorphism classes of right ideals. The group $G = \mr{Aut}_{\C}(A_1)$ acts on the set $\mc{C}$. A key observation by Y. Berest and G. Wilson, which was the foundation of their seminal work \cite{BerestWilson1}, is that this set is precisely the ``adelic Grassmanian" appearing in integrable systems. They went on to show that  
\[
\mc{C} = \bigsqcup_{n \ge 0} \mc{C}_n,
\]
where each $\mc{C}_n$ is a $G$-orbit which can also be identified with Wilson's compactified Calogero-Moser space \cite{Wilson}. This connection between right ideals in the Weyl algebra and the Calogero-Moser space is very mysterious and several different constructions of the above bijection have appeared in the literature in order to try and demystify the connection. Since each $\mc{C}_n$ is a $G$-orbit, Corollary~\ref{cor:Moritaequiv} says that the ideals in $\mc{C}_n$ all give rise to isomorphic endomorphism rings. Thus, right ideals in the Weyl algebra provide a countable family $E_m$, $m \in \N$, of rings, all Morita equivalent but pairwise non-isomorphic. Stafford shows that $E_m$ can be identified with the ring of differential operators on the non-normal rational curve $X_m$ where $\mathcal{O}(X_m) = \C \oplus x^{m+1} \C[x]$. 

The group $G = \mr{Aut}_{\C}(A_1)$ was first studied by Dixmier \cite{DixmierWeyl}, who gave an explicit set of generators. %; see also \cite[\S11]{BerestWilson1} or \cite{BEETransitive}. 
Setting $G_m = \mr{Aut}_{\C}(E_m)$, Stafford asked \cite[p. 636]{StaffordEdnomorphisms}: 
\begin{enumerate}
\item Does there exist a description of $G_m$ for $m >0$ similar to Dixmier's for $G$?
\item Is it true that $G_m \not\cong G$ when $m > 0$?
\end{enumerate}
Berest--Eshmatov--Eshmatov \cite{BEEDixmier} show that the answer to both questions is ``yes".  For higher $n$, the group $\mr{Aut}_{\C}(A_n)$ is not so well-understood. Belov-Kanel and Kontsevich conjectured in \cite{KontAutoConj} that  $\mr{Aut}_{\C}(A_n)$ should be equal to the group of Poisson $\C$-algebra automorphisms of the polynomial ring $P$. This conjecture has been proven in \cite{BelovKanelProof}, and a related (but slightly weaker) conjecture by Belov-Kanel and Kontsevich \cite[Conjecture~6]{KontAutoConj}  was proven recently by C. Dodd \cite{Doddpcycle}.  

%Closely related to these questions is Dixmier's (famous) conjecture \cite{DixmierWeyl} which states that every endomorphism of $A_n$ is an automorphism, $\mr{End}_{\C}(A_n) = \mr{Aut}_{\C}(A_n)$. For $n = 1$, the conjecture  was recently proven by A. Zheglov \cite{zheglov2025conjecturedixmierweylalgebra}. Dixmier's conjecture is known to imply the infamous Jacobian conjecture; see for instance \cite[Theorem~4.2]{PrimerDmodules}.

\section{Generalizations to quantized symplectic singularities}

In this final section we wish to advertise the fact that Stafford's results, and the subsequent ones by Bernstein-Lunts, should have natural generalizations to quantized symplectic singularities. These algebras are currently the subject of intense research within geometric representation theory. 

A normal affine variety $X$ is said to have symplectic singularities if its smooth locus has a holomorphic symplectic form compatible, in a strong sense, with all resolutions of singularities. See the original paper \cite{Beauville} by Beauville for details. We say that a variety $X$ with symplectic singularities is \textit{conic} if $\C[X]$ is connected graded such that the symplectic form is homogeneous of (strictly) positive weight. Let $U = \bigcup_{n \ge 0} \mc{F}_n U$ be a filtered $\C$-algebra such that $\gr_{\mc{F}} U$ is a commutative ring. Then $\gr_{\mc{F}} U$ is a Poisson algebra. If $\gr_{\mc{F}} U \cong \C[X]$ as graded Poisson algebras, where $X$ is a conic symplectic singularity, then $U$ is said to be a quantization of a symplectic singularity. 

The quantizations of a given conic symplectic singularity $X$ have been classified by I. Losev \cite[Theorem~3.4]{LosevDeformations}. In particular, the Weyl algebra is the unique quantization of $\C^{2n}$. The variety $X$ has a finite stratification by symplectic leaves - these are smooth locally closed Poisson subvarieties where the Poisson structure is non-degenerate. If $M$ is a finitely generated $U$-module then the characteristic variety $\Ch(M)$ of $M$ is coisotropic by Gabber's Theorem. This means that $\Ch (M) \cap \mc{L}$ is a coisotropic subvariety for all leaves $\mc{L}$. We say that $M$ is \textit{holonomic} if $\Ch (M) \cap \mc{L}$ is a Lagrangian subvariety for every leaf $\mc{L}$. 

%\subsection{Bernstein's inequality revisited} In \cite{BernsteinLosev}, Losev gives a vast generalization of the Bernstein inequality\footnote{At the time of writing, there is a gap in the proof of Lemma~3.3 in \cite{BernsteinLosev}. This lemma is a key step in the proof of the main result \cite[Theorem~1.1]{BernsteinLosev}.} to any quantization $U$ of a symplectic singularity (in fact, he works in a more general setting). 

%\begin{thm}\cite[Theorem~1.1]{BernsteinLosev}\label{thm:LosevB1}
%	Suppose that $M$ is a simple $U$-module with annihilator $I$.
%	\begin{enumerate}
%		\item $2 \dim \Ch(M) \ge \dim \Ch(U/I)$. 
%		\item If $\mc{L}$ is a symplectic leaf such that $\Ch(M) \cap \mc{L} \subset \Ch(M)$ is open non-empty, then $\Ch(U/I) = \overline{\mc{L}}$ and the intersection $\Ch(M) \cap \mc{L}$ is a coisotropic subvariety in $\mc{L}$.
%\item Suppose that the $U$-bimodule $U$ has finite length. Then $2 \dim \Ch(M) \ge \dim \Ch(U/I)$.
%\item The intersection of any irreducible component of $\Ch(M)$ with $\mc{L}$ is non-empty.
%\end{enumerate}
%\end{thm}

%In the case of holonomic modules, the above theorem can be strengthen to give:

%\begin{thm}\cite[Theorem~1.2]{BernsteinLosev}
 %Suppose that $N$ is a holonomic $U$-module with annihilator $I$. 
%		\begin{enumerate}
%\item If $N$ is simple, then $2 \dim \Ch(N) = \dim \Ch(U/I)$.
%\item If the $U$-bimodule $U$ has finite length, then $N$ has finite length. Consequently, 
%\[
%2 \dim \Ch(N) = \dim \Ch(U/I).
%\]
%\item If $N$ is simple, then the variety $\Ch(N)$ is equi-dimensional.
%\end{enumerate}
%\end{thm}

If $M$ is a simple $U$-module then, by definition, its annihilator is a primitive ideal $I \lhd U$. V.~Ginzburg \cite[Theorem~2.1]{Primitive} showed that the associated variety $\Ch(U/I) \subset X$ is the closure of a symplectic leaf. Notice that $\Ch(M) \subset \Ch(U/I)$. We ask:

\begin{question}\label{que:q1}
If $\Ch(U/I)$ is the closure of the leaf $\mc{L}$ then is $\mc{L} \cap \Ch(M) \neq \emptyset$? 
\end{question}

If $\Ch(M)$ intersects $\mc{L}$ then $\mc{L} \cap \Ch(M)$ is a coistropic subvariety of $\mc{L}$ and hence $2 \dim \Ch(M) \ge \dim \mc{L}$. This implies that $2 \GKdim(M) \ge \GKdim(U/I)$ i.e., the Bernstein inequality holds. If $\Ch(U/I) = \overline{\mc{L}}$ then we say that $I$ is a quantization of $\mc{L}$ over $U$. For a given $U$, not all leaves of $X$ admit a quantization. In fact, for a given leaf $\mc{L}$ there might not be any quantization $U$ of $X$ such that $\mc{L}$ can be quantized over $U$. 

\begin{question}
Which leaves of $X$ can be quantized over $U$? 
\end{question}

Comparing Question~\ref{que:q1} with the results on existence of non-holonomic simple modules for Weyl algebras, it is natural to ask:

\begin{question}
Assume that the leaf $\mc{L}$ can be quantized over $U$. Is it the case that for all $0 \le i \le (1/2) \dim \mc{L}$ there exists a simple $U$-module $M_i$ such that 
	\[
	\Ch(U / \ann_U M_i ) = \overline{\mc{L}} \quad \textrm{ and } \quad \GKdim(M_i) = (1/2) \dim \mc{L} + i ?
	\] 
\end{question}

A related question concerns the behaviour of Ext-groups between simple modules. 

\begin{question}
Let $M,N$ be simple $U$-modules. 
	\begin{enumerate}
	\item Is there a concise characterization of when the group $\Ext^i_U(M,N)$ is a finite-dimensional $\C$-vector space?  
	\item In particular, if $M,N$ are holonomic is it the case that $\Ext^i_U(M,N)$ is always a finite-dimensional $\C$-vector space?
	\end{enumerate}
\end{question}

\subsection{Simplicity} We have seen throughout that a key property of the Weyl algebra is that it is a simple ring. Quantizations $U$ are not always simple rings. However, in Losev's classification \cite{LosevDeformations} there is always a preferred quantization of $X$ called the canonical quantization, see \cite[\S5]{losev4unipotentidealsharishchandrabimodules}.  

\begin{conjecture}\cite[Conjecture~5.5.1]{losev4unipotentidealsharishchandrabimodules}\label{conj:losevsimple}
	The canonical quantization is simple. 
\end{conjecture}

If $X$ has $\Q$-factorial terminal singularities then there is a unique quantization $U$ of $X$. In this case, the above conjecture says that $U$ should be simple. This is a special case of an earlier conjecture by Losev \cite[Conjecture~3.1]{LosevDerivedSRA}. The truth of Conjecture~\ref{conj:losevsimple} would imply that every conic symplectic singularity has a simple quantization. If $U$ is a simple quantization then it would be interesting to know if the results of \cite{StaffordModuleStructure} can be shown to hold for $U$. 

\subsection{Automorphisms} In Section~\ref{sec:subsequent}, we saw the appearance of the automorphism group of the Weyl algebra in relation to right ideals in the first Weyl algebra. Therefore, it is natural to consider the group of $\C$-algebra automorphisms $\mr{Aut}_{\C}(U)$ of a quantization $U$ of a conic symplectic singularity $X$. Since Losev's classification says that there is actually a family of quantizations, parametrized by a certain space $\mathfrak{h}^* / W$ where $W$ is Namikawa's Weyl group, it is natural to consider instead the groupoid $\mathsf{QIso}(X)$ whose objects are points of $\mathfrak{h}^* / W$ and, for $c,c' \in \mathfrak{h}^* / W$, $\Hom_{\mathsf{QIso}(X)}(c,c')$ is the set of all isomorphisms $U_c \iso U_{c'}$ (not necessarily filtration preserving). As explained in \cite[\S2]{LosevDeformations}, the same space $\mathfrak{h}^* / W$ also parametrizes all filtered Poisson deformations of $\C[X]$ and one can similarly form the groupoid $\mathsf{PIso}(X)$ whose objects are points of $\mathfrak{h}^* / W$ and $\Hom_{\mathsf{PIso}(X)}(c,c')$ is the set of all Poisson isomorphisms $A_c \iso A_{c'}$ between the filtered Poisson deformations $A_c,A_{c'}$ (again, not necessarily filtration preserving). The following question was raised in \cite[Final~Remarks]{Castellan}. 

\begin{question}
Are the groupoids $\mathsf{QIso}(X)$ and $\mathsf{PIso}(X)$ isomorphic over $\mathfrak{h}^* / W$? 
\end{question}

S. Castellan \cite{Castellan} shows that the answer is yes when $X$ is a Kleinian singularity of type $A$ or $D$ and conjectures it should be the case for any Kleinian singularity. The equivalence $\mathsf{QIso}(\C^{2n}) \simeq \mathsf{PIso}(\C^{2n})$ is precisely Belov-Kanel and Kontsevich's conjecture since $\h^*/W = \{ 0 \}$ in this case. 

%\begin{remark}
%	A "Conze embedding" can be shown to exist for many quantizations of conic symplectic singularities so, in principal, one should be able to use Stafford's results to prove the existence of non-holonomic simple modules for a large class of quantizations just as was done for the primitive quotient of $\mf{sl}_2 \times \mf{sl}_2$ in Section~\ref{sec:applicationstoenvelopin}. 
%\end{remark}

%\bibliographystyle{plain}
%\bibliography{biblo}

\def\cprime{$'$} \def\cprime{$'$} \def\cprime{$'$} \def\cprime{$'$}
  \def\cprime{$'$} \def\cprime{$'$} \def\cprime{$'$} \def\cprime{$'$}
  \def\cprime{$'$} \def\cprime{$'$} \def\cprime{$'$} \def\cprime{$'$}
  \def\cprime{$'$} \def\cprime{$'$}

\end{document}